\documentclass[12pt]{amsart}
\usepackage{amsmath,amsthm,amssymb}
\usepackage[matrix,arrow,curve]{xy}
\makeatletter\@addtoreset{equation}{section} \makeatother

\newtheorem{theorem}[equation]{Theorem}

\newtheorem{statement}[equation]{Statement}
\newtheorem{lemma}[equation]{Lemma}
\newtheorem{corollary}[equation]{Corollary}
\newtheorem{conjecture}[equation]{Conjecture}

\theoremstyle{definition}
\newtheorem{example}[equation]{Example}

\theoremstyle{remark}
\newtheorem{remark}[equation]{Remark}

\def \Q {\mathbb{Q}}
\def \P {\mathbb{P}}
\def \Z {\mathbb{Z}}
\def \C {\mathbb{C}}
\def \F {\mathbb{F}}
\def \Y {\mathcal{Y}}
\def \H {\mathcal{H}}
\def \O {\mathcal{O}}
\def \QQ {\mathcal{Q}}

\def \rk {\mathrm{rk}}
\def \Bir {\mathrm{Bir}}
\def \Aut {\mathrm{Aut}}
\def \PGL {\mathrm{PGL}}
\def \Pic {\mathrm{Pic}}

\def \ge {\geqslant}
\def \le {\leqslant}

\newcommand{\tit}{
$\Q$-factorial quartic threefolds
}
\author{C.\,Shramov}
\title{\tit}

\begin{document}

\maketitle

\begin{abstract}
We prove that a nodal quartic threefold $X$ containing no planes is 
$\Q$-factorial provided that it has not more than $12$ singular points,
with the exception of a quartic with exactly $12$ singularities containing 
a quadric surface. We give some geometrical constructions related to the latter
quartic.
\end{abstract}

\section{Introduction}

A quartic threefold is a classical object of interest in birational geometry. 
The group of birational automorphisms of a nonsingular quartic was first
studied in~\cite{IskovskikhManin}; birational automorphism groups of
singular quartics were studied in particular in~\cite{Pukhlikov-quartic} 
(the case of a single simple double point) and~\cite{Mella} 
(the case of arbitrary number of simple double points). M.\,Mella 
(see~\cite{Mella}) proved birational rigidity and described the generators
of the group of birational automorphisms of any $\Q$-factorial nodal
quartic. On the other hand, one cannot extend these results to 
the case of non-$\Q$-factorial varieties: 
there are non-$\Q$-factorial nodal quartics threefolds that are rational.
This motivates the following question: which conditions guarantee 
$\Q$-factoriality of a nodal quartic?

Recall that a variety $V$ is called \emph{factorial} if any Weil 
divisor on $V$
is also a Cartier divisor, 
and \emph{$\Q$-factorial} if an appropriate multiple of any Weil divisor 
is a Cartier divisor. A singular point $x\in V$ is called 
\emph{a simple double point} (or \emph{a simple double singularity}) 
if it has a neighborhood analytically isomorphic to a neighborhood
of a vertex of a cone over a nonsingular quadric. If $V$ is a hypersurface in 
$\P^n$ given by the equation $F=0$, the latter is equivalent to non-degeneracy
of the Hessian $H(F)$ at $x$. The variety with only simple double
singularities is called \emph{nodal}.

Let $X$ be a nodal quartic threefold. Then the following is well-known.

\begin{theorem}[see~\cite{Cheltsov-quartic}]\label{theorem:no-plane}
If the number of singular points of $X$ does not exceed $8$, then 
$X$ is $\Q$-factorial; 
if it is equal to $9$, and $X$ contains no planes, then $X$ 
is again $\Q$-factorial.
\end{theorem}

This result is a particular case of the following conjecture 
(see, for example,~\cite{Ciliberto}, \cite{Cheltsov-quartic}, 
and also~\cite{Cheltsov-points} for some advances).

\begin{conjecture}\label{conjecture:CdG}
Let $Y$ be a nodal hypersurface of degree $d$ in $\P^4$ with\nolinebreak $s$ 
singular points. If $s<(d-1)^2$, then $Y$ is $\Q$-factorial; if $Y$ contains 
no planes and  $s<2(d-1)(d-2)$, then $Y$ is $\Q$-factorial; if 
$s\nolinebreak\le\nolinebreak 2(d-\nolinebreak 1)\nolinebreak 
(d-\nolinebreak 2)$, 
and $Y$ contains neither planes nor quadrics, then $Y$ 
is again $\Q$-factorial.
\end{conjecture}

By now Conjecture~\ref{conjecture:CdG} is checked only for $d\le 3$.

The main goal of this paper is to show the following.

\begin{theorem}\label{maintheorem}
Let a nodal quartic threefold $X$ with not more than $11$ singularities 
contain no planes. Then $X$ is $\Q$-factorial. If $X$ has $12$ singularities,
then $X$ is again $\Q$-factorial with the exception of the case when 
$X$ contains a quadric surface.
\end{theorem}

Theorem~\ref{maintheorem} generalizes Theorem~\ref{theorem:no-plane};
it is a particular case of Conjecture~\ref{conjecture:CdG}. Since there 
are non-$\Q$-factorial nodal quartics containing no planes 
and having $12$ singular points (the example is given by
a general quartic containing a quadric --- 
see Example~\ref{example:construction} below), 
the statement of Theorem~\ref{maintheorem} cannot be refined (as well as
the general statement of Conjecture~\ref{conjecture:CdG}).

The proof of Theorem~\ref{maintheorem} occupies 
section~\ref{section:Q-factoriality}; it uses some simple consequences
of Bese theorem (see Theorem~\ref{theorem:Bese}) that are collected 
in section~\ref{section:auxiliary-surfaces}. 
In section~\ref{section:constructions}
we consider some constructions related to a general quartic containing
a quadric: we describe a (well-known) transformation of this quartic
to a complete intersection of a quadric and a cubic in $\P^5$ with a single 
double point (Example~\ref{example:construction}); 
we also give an example of a birational involution of the 
latter variety that does not exist in the nonsingular case 
(Example\ref{example:involution} and Statement~\ref{statement:action}). 

All the varieties are defined over the field of complex numbers $\C$.

\smallskip
The author is grateful to I.\,Cheltsov and S.\,Galkin for useful discussions.

%%%%%%%%%%%%%%%%%%%%%%%%%%%%%%%%%%%%%%%%

\section{Auxiliary statements about points on surfaces}
\label{section:auxiliary-surfaces}

In this section we collect some statements that will be used
in section~\ref{section:Q-factoriality}.

\begin{theorem}[see~\cite{Bese}]\label{theorem:Bese}
Let $\pi:\F_{r, N}\to\F_r$ be a blow-up of a rational ruled surface
$\F_r$ at $N$ points $p_1, \ldots, p_N$. 
Let $D\sim as+bf$ be a divisor on $\F_r$,
where $s$ and $f$ denote the classes of an exceptional section and a fiber
of $\F_r$ respectively, and $a, b>0$. Let $L_N=\pi^*D-\sum_{i=1}^N E_i$, where 
$E_i$ are exceptional divisors of $\pi$.  Let $\rho=\dim |D-K_{\F_r}|$, and 
$$h=h^1(\F_r, \O_{\F_r}(D)\otimes K_{\F_r}^{-1})=
\rho-\frac{D^2-3DK_{\F_r}+16}{2}.$$ 
Then the line bundle $\O_{\F_{r, N}}(L_N)$ is generated by its global sections, provided that the following conditions hold:

(i) $N\le \lfloor\frac{\rho-4}{3}\rfloor$,

(ii) $D^2\ge 7+4h$,

(iii) any curve $C$ of bidegree $(x, y)\neq (0, 0)$ (i.\,e. linearly 
equivalent to $xs+yf$), $0\le x\le a+2$,
$0\le y\le \lfloor\frac{b+2+r}{2}\rfloor$, contains not more than 
$C(D-K_{\F_r}-C)-2$ points of the set $\{p_1, \ldots, p_N\}$.
\end{theorem}

\begin{corollary}\label{corollary:10-on-p1p1}
Let $p_1, \ldots, p_{10}, q\in \F_0=\P^1\times\P^1$ be such points that 
no curve of bidegree $(1, 0)$ or $(0, 1)$ contains $4$ points of 
$p_1, \ldots, p_{10}$, no curve of bidegree $(1, 1)$ contains $7$ 
of them and no curve of bidegree $(2, 1)$ or $(1, 2)$ containes all of them. 
Then there is a divisor $D$ of bidegree $(3, 3)$ passing through
$p_1, \ldots, p_{10}$ and not passing through\nolinebreak $q$. 
\end{corollary}
\begin{proof}
In the notations of Theorem~\ref{theorem:Bese} we have 
$\rho=35$, $h=0$, $D^2=18$,
hence the conditions (i) and (ii) hold. Condition (iii) claims that
any curve  
$C$ of bidegree $(x, y)\neq (0, 0)$, $0\le x\le 5$, $0\le y\le 2$,
contains not more than  
$$C(D-K_{\P^1\times\P^1}-C)-2=5x+5y-2xy-2$$
points of $p_1, \ldots, p_{10}$; the check of the latter is straightforward 
for all the pairs $(x, y)$.
\end{proof}

\begin{corollary}\label{corollary:10-on-f2}
Let $p_1, \ldots, p_{10}, q\in \F_2$ be such points that 
the exceptional section contains no points of  
$p_1, \ldots, p_{10}$, no fiber contains $4$ of them, no curve
of bidegree $(1, 2)$ contains $7$ of them and no curve of bidegree $(1, 3)$ 
containes all of them. 
Then there is a divisor $D$ of bidegree $(3, 6)$ passing through 
$p_1, \ldots, p_{10}$ and not passing through $q$. 
\end{corollary}
\begin{proof}
In the notations of Theorem~\ref{theorem:Bese} we have 
$\rho=35$, $h=0$, $D^2=18$,
hence the conditions (i) and (ii) hold. Condition (iii) claims that 
any curve $C$ of bidegree $(x, y)\neq (0, 0)$, $0\le x\le 5$, $0\le y\le 5$, 
contains not more than $2x^2+5y-2xy-2$ points of $p_1, \ldots, p_{10}$ and is
a straightforward consequence of the assumptions (note that condition (iii) for 
$(x, y)=(1, 0)$ means that the exceptional section does not contain any of the
points $p_i$).
\end{proof}

\begin{corollary}\label{corollary:8-on-f2}
Let $p_1, \ldots, p_8, q\in \F_2$ be such points that
the exceptional section contains no points of 
$p_1, \ldots, p_8$, no fiber contains $4$ of them, 
no curve of bidegree $(1, 2)$ contains $6$ of them and no curve of bidegree
$(1, 3)$ containes all of them. 
Then there is a divisor $D$ of bidegree $(2, 5)$ passing through 
$p_1, \ldots, p_8$ and not passing through $q$. 
\end{corollary}
\begin{proof}
In the notations of Theorem~\ref{theorem:Bese} we have 
$\rho=29$, $h=0$, $D^2=12$,
hence conditions (i) and (ii) hold. Condition (iii) claims that
any curve $C$ of bidegree $(x, y)\neq (0, 0)$, $0\le x\le 4$, $0\le y\le 4$, 
contains not more than $2x^2+x+4y-2xy-2$ points of $p_1, \ldots, p_8$
and is a straightforward consequence of the assumptions.
\end{proof}

\begin{corollary}\label{corollary:10-on-a-quadric-in-p3}
Let $Y\subset\P^3$ be an irreducible quadric; let 
$p_1, \ldots, p_{10}, q\in Y$ be such points that 
no line contains $4$ of the points $p_1, \ldots, p_{10}$ 
no conic contains $7$ of them and no twisted cubic containes all of them.
Then there is a divisor $D\in\left.\O_{\P^3}(3)\right|_{Y}$ passing 
through $p_1, \ldots, p_{10}$ and not passing through $q$.
\end{corollary}
\begin{proof}
If $Y$ is a nonsingular quadric, then the statement 
follows from Corollary~\ref{corollary:10-on-p1p1}. If $Y$ is 
a cone over a conic and its vertex does not coincide with any of the points
$p_i, q$, then we are done by Corollary~\ref{corollary:10-on-f2} applied
to the preimages of $p_i, q$ under the resolution map $g:\F_2\to Y$. 
If $q$ is the vertex, we apply Corollary~\ref{corollary:10-on-f2} to the
points $g^{-1}(p_i)$ and an arbitrary point $\tilde{q}\in s\subset\F_2$; 
this gives us a divisor $\tilde{D}\sim 3s+6f$ passing through
$g^{-1}(p_i)$ and not containing the exceptional section $s$, 
and hence not intersecting $s$, since $\tilde{D}s=0$. The divisor 
$D=g_*\tilde{D}$ is just what we looked for.

Now we have to consider the case when one of the points $p_i$ (say, $p_{10}$)
coincides with the vertex of $Y$. Since not all of the points $p_i$ lie on 
a single line, we may choose a line $l\subset Y$ such that $q\notin l$
and one of the points 
$p_1, \ldots, p_9$ is contained in $l$. Moreover, if there is a conic
containing $6$ of the points $p_1, \ldots, p_9$, then such a conic is unique,
and we can choose $l$ so that it contains one of these $6$ points 
(say, $p_9$ again). 
So we may assume that no $3$ of the points of $p_1, \ldots, p_8$ 
lie on a line, no $6$ of them lie on a conic and all $8$ are not contained
in a twisted cubic.
Applying Corollary~\ref{corollary:8-on-f2} to the points 
$g^{-1}(p_1), \ldots,
g^{-1}(p_8), g^{-1}(q)$, we find a divisor $\tilde{D}\sim 2s+5f$ passing 
through $g^{-1}(p_1), \ldots, g^{-1}(p_8)$ and not passing through  
$g^{-1}(q)$. Now we construct the desired divisor as $D=g_*\tilde{D}\cup l$.
\end{proof}

%%%%%%%%%%%%%%%%%%%%%%%%%%%%%%%%%%%%%%%%%%%%%%%%%%%%%%%%%%%%%%%%%%%%%%%
\section{$\Q$-factoriality}\label{section:Q-factoriality}

Throughout this section $X$ denotes a nodal quartic threefold
containing no plains and singular at the points 
$p_1, \ldots, p_s$.

In the case of a nodal threefold $\Q$-factoriality is equivalent to 
factoriality; on the other hand, factoriality of a nodal Fano threefold
$X$ is equivalent to a topological condition $\rk H^2(X, \Z)=\rk H_4(X, \Z)$ 
(in our case $\rk H^2(X, \Z)=\rho(X)=1$).
Hence to prove $\Q$-factoriality of $X$ it suffices to check that 
for a small resolution $h:\tilde{X}\to X$ we have
$\rho(\tilde{X})=s+1$. To check this it is sufficient
to show that, using the terminology of~\cite{Cynk}, 
the defect of $X$ equals zero (cf. the proof of Theorem~2  
in~\cite{Cynk}, and also~\cite{CheltsovPark}). The latter condition
means that the singular points $p_1, \ldots, p_s$ impose independent
conditions on the hypersurfaces of degree $3$ in $\P^4$, i.\,e. that 
for any point $p_i$, $1\le i\le s$, there is a cubic hypersurface 
$D_i\subset\P^4$ such that $p_i\not\in D_i$ and $p_j\in D_i$ for $j\neq i$. 

The main tool for checking this independency condition is provided 
by the following theorem proved in~\cite{EisenbudKoh}.

\begin{theorem}\label{theorem:Eisenbud}
The points $p_1, \ldots, p_s\in\P^n$ impose independent conditions on forms 
of degree $d$ if each linear subspace of dimension $k$ contains 
at most $dk+1$ of the points $p_1, \ldots, p_s$.
\end{theorem}

The points $p_1, \ldots, p_s$ have the following property.

\begin{lemma}\label{lemma:6-points}
No line contains $4$ of the points $p_1, \ldots, p_s$,
no plane contains $7$, and no twisted cubic contains $10$ 
of them. If the points  
$p_1, \ldots, p_6$ are contained in a plane $P$, then either $X$ intersects 
$P$ along a (possibly reducible) double conic, or $X$ intersects
$P$ along a union of four different lines $l_1, \ldots, l_4$ and the points
$p_i$, $1\le i\le 6$, coincide with the points of pairwise intersections of the 
lines $l_j$, $l_k$. 
\end{lemma}
\begin{proof}
The first three statements hold regardless to the assumption that 
$X$ contains no planes (see, for example,~\cite[Lemma~2.7]{Cheltsov-points}). 
Let the points $p_1, \ldots, p_6$ lie on a plane $P$.
Assume that $p_1, \ldots, p_6$ lie on a conic $Q$.
If there exists a point $q\in (X\cap P)\setminus Q$, then there is a 
two-dimensional family of plane cubics passing through $p_1, \ldots, p_6, q$.
Let $C$ be a general curve of this family. Then $C$ is irreducible 
and $C\subset X$, since otherwise the intersection index of $C$ and $X$ in 
$\P^4$ would be at least $13$.
As $C$ is movable, we have $P\subset X$ that is a contradiction.

Now assume that $X\cap P$ is neither a conic nor a union of four lines. 
It means that there is a point $q\in X\cap P$ such that 
the points $p_1, \ldots, p_6, q$ are not contained in a single conic, 
and no $4$ of them are collinear. So there is again a 
two-dimensional family of plane cubics passing through $p_1, \ldots, p_6, q$;
a general curve of this family is irreducible and hence is contained in
$X$, a contradiction.
\end{proof}

Theorem~\ref{theorem:Eisenbud} and Lemma~\ref{lemma:6-points} 
immediately imply

\begin{corollary}\label{corollary:10-points}
Let $s\le 10$. Then $X$ is $\Q$-factorial.
\end{corollary}

On the other hand, Lemma~\ref{lemma:6-points} gives the following 

\begin{corollary}\label{corollary:reducible-quadric}
Assume that $s\ge 11$ and the points $p_1, \ldots, p_s$
are contained in a reducible quadric surface (i.\,e. in a pair of
planes which span a three-dimensional subspace). Then either 
$p_1, \ldots, p_s$ impose independent conditions on the forms of degree $3$,
or $p_1, \ldots, p_s$ are also contained in an irreducible quadric surface.
\end{corollary}
\begin{proof}
Let $p_1, \ldots, p_s$ be contained in a union of the planes $L$ and $L'$ 
so that $p_1, \ldots, p_6\in L$, $p_7, \ldots, p_s\in L'$,
$L\cap L'=l$. Note that $s\le 12$, and at most one of the points $p_i$ 
lie on the line $l$.
In particular, $p_1, \ldots, p_6$ are contained in a conic if and only if
$p_7, \ldots, p_s$ are, since our assumptions 
imply that the union of four lines that intersect in singular points of $X$
intersects the line $l$ at least by three points, and a conic intersects $l$ 
at most by two points.

Assume that the points $p_1, \ldots, p_6$ lie on a conic $Q_1$, and
the points $p_7, \ldots, p_s$ lie on a conic $Q_2$.
Note that $Q_1\cap l=Q_2\cap l$ (otherwise $X$ would contain one of the planes 
$L$ or $L'$), i.\,e. the conics $Q_i$ intersect $l$ by the same pair
of (possibly coinciding) points. Any two conics with such a property 
are contained in an irreducible quadric.

Hence if the points $p_1, \ldots, p_s$ are not contained in an irreducible
quadric, Lemma~\ref{lemma:6-points} implies that the intersection of $X$ 
and $L$ is a union of four lines $l_1, \ldots, l_4$.
So the intersection of $X$ and $L'$ is also a union of four lines;
denote them by $l_1', \ldots, l_4'$. We may assume (renumbering the 
lines $l_i'$ if necessary) that $l_i\cap l=l_i'\cap l$. Now it is easy to
see that the points $p_1, \ldots, p_s$ impose independent conditions
on the forms of degree $3$.
(For example, to construct a cubic surface $S_1$ passing through 
$p_2, \ldots, p_s$ and not passing through $p_1$, we may proceed as follows:
if $p_1$ is not contained in $l_1$ and $l_2$, then $S_1$ is a union
of the planes $P_1=\langle l_1, l_1'\rangle$, $P_2=\langle l_2, l_2'\rangle$
and some plane $P_3$, passing through the only point of 
$p_7, \ldots, p_s$ that lies on neither $l_1'$ nor $l_2'$, and
passing through no other singular point.)
\end{proof}

We'll also need the following observations. 

\begin{lemma}\label{lemma:no-quadric}
Assume that $s\le 12$ and the points $p_1, \ldots, p_s$ 
are contained in a three-dimensional subspace $\P^3\subset\P^4$
but not in a quadric surface. Then $p_1, \ldots, p_s$ impose 
independent conditions on the forms of degree $3$ in $\P^3$ 
(and hence in $\P^4$ as well).
\end{lemma}
\begin{proof}
We may assume that $s=12$. It suffices to find a cubic surface,
passing through $p_1, \ldots, p_{11}$ and not passing through $p_{12}$. 
Since the points $p_1, \ldots, p_{12}$ are not contained in a single 
quadric, we may find a quadric containing $9$ of them 
(say, $p_1, \ldots, p_9$) and not containing $p_{12}$. If the points 
$p_{10}$, $p_{11}$ and $p_{12}$ are noncollinear,
we may construct the desired cubic surface 
as a union of this quadric with a plane
passing through $p_{10}$ and $p_{11}$ and not passing through $p_{12}$. 
If the points $p_{10}$, $p_{11}$ and 
$p_{12}$ are collinear, Lemma~\ref{lemma:6-points} implies that the points 
$p_i$, $p_{11}$ and $p_{12}$ are not collinear for any $1\le i\le 9$. 
On the other hand, not all the quadrics passing through the collections 
of the points 
$\{p_1, \ldots, p_9,\nolinebreak p_{10}\nolinebreak\}
\setminus\nolinebreak\{p_i\}$, $1\le i\le 9$, also 
pass through $p_{12}$, and so we may again obtain our cubic 
as a union of a quadric and a plane.
\end{proof}

\begin{lemma}
\label{lemma:11-on-quadric}
Assume that $s\le 11$ 
and the points $p_1, \ldots, p_s$ are contained in a quadric surface.
Then $p_1, \ldots, p_s$ impose independent conditions on the forms of degree
$3$ in $\P^4$.
\end{lemma}
\begin{proof}
Due to Corollary~\ref{corollary:reducible-quadric} we may assume that 
the points $p_1, \ldots, p_s$ are contained in an irreducible quadric. 
Corollary~\ref{corollary:10-on-a-quadric-in-p3} finishes the proof.
\end{proof}

\begin{corollary}\label{corollary:11-and-1}
Assume that $s=12$ and the points $p_1, \ldots, p_{11}$ span a three-dimensional
subspace not containing $p_{12}$. Then $p_1, \ldots, p_{12}$ impose independent 
conditions on the forms of degree $3$ in $\P^4$. 
\end{corollary}
\begin{proof}
Apply Lemma~\ref{lemma:no-quadric} and Lemma~\ref{lemma:11-on-quadric}.
\end{proof}

\begin{lemma}
\label{lemma:12-on-quadric}
Assume that $s=12$ and the points $p_1, \ldots, p_{12}$ are contained
in a quadric surface.
Then either $p_1, \ldots, p_{12}$ impose independent conditions on 
the forms of degree $3$ in $\P^4$, or $p_1, \ldots, p_{12}$ are 
contained in a pencil of quadric surfaces in some three-dimensional
subspace.
\end{lemma}
\begin{proof}
Due to Corollary~\ref{corollary:reducible-quadric} we may assume that 
the points $p_1, \ldots, p_s$ are contained in an irreducible 
quadric $Q$; we may also assume that $s=12$. Now assume 
that any quadric passing through $p_1,\ldots, p_{12}$ coincides with $Q$
(otherwise the second assertion holds). 
Let us prove that there exists a (reducible) cubic surface 
passing through 
$p_1, \ldots, p_{11}$ and not passing through $p_{12}$. Let $Q'$ be
a quadric different from $Q$, such that $Q'$ contains $8$ of the points 
$p_i$ (say, $p_1, \ldots, p_8$) and does not contain $p_{12}$. If the points
$p_9, \ldots, p_{12}$ are noncoplanar, we may construct the 
cubic surface in question as a union of a quadric $Q'$ and a plane passing 
through $p_9$, $p_{10}$ and $p_{11}$. So let the points 
$p_9, \ldots, p_{12}$ lie in a plane $P$; renumbering
the points $p_9$, $p_{10}$, $p_{11}$ if necessary, we may assume that 
$P$ is spanned by $p_{10}$, $p_{11}$ and  
$p_{12}$. By Lemma~\ref{lemma:6-points} this plane contains at most
two of the points $p_1, \ldots, p_8$. Let these be $p_7$ and $p_8$ 
(a similar argument is applicable to the cases when $P$ contains 
either only one of the points 
$p_1, \ldots, p_8$ or none of them). 
If all the quadrics passing through the collections 
$\{p_1, \ldots, p_8, p_9\}\setminus\{p_i\}$, $1\le i\le 6$, 
also pass through the point $p_{12}$, then the following condition holds: 
any quadric, containing the points $p_7$, $p_8$ and $p_9$ also
contains $p_{12}$,
i.\,e. the points $p_7$, $p_8$, $p_9$ and $p_{12}$ are collinear;
the latter contradicts Lemma~\ref{lemma:6-points}. Otherwise there exists,
say, a quadric $Q_1$ containing $p_2, \ldots, p_9$ and not containing
$p_{12}$; now we obtain a desired cubic as a union of $Q_1$ with a plane 
passing through $p_1$, $p_{10}$ and $p_{11}$.
\end{proof}

\begin{lemma}\label{lemma:contains-quadric}
Assume that $s=12$ and the points $p_1, \ldots, p_{12}$ are contained
in a pencil of quadric surfaces in some three-dimensional subspace. 
Then $X$ contains a quadric surface.
\end{lemma}
\begin{proof}
Let $Q$ and $Q'$ be two different quadric surfaces containing 
$p_1, \ldots, p_{12}$; let a curve $E$ be the intersection of $Q$ and $Q'$. 
We may assume that $Q$ and $Q'$ are irreducible and $Q$ contains 
a point $q\in X$, $q\not\in E$. Let $\H$ be a linear system
on $Q$ that consists of the restrictions of the cubic forms vanishing at
the points $p_1, \ldots, p_{12}$ and $q$. Lemma~\ref{lemma:6-points} 
implies that $\H$ has no base components. But the intersection index of
the corresponding curves with $X$ in $\P^4$ is at least $25$, hence 
$Q\subset X$.
\end{proof}

\begin{corollary}\label{corollary:12-points}
Assume that $s=12$ and $11$ of the points $p_1, \ldots, p_{12}$ are contained
in a three-dimensional subspace. Then $X$ is either 
$\Q$-\nolinebreak factorial, 
or contains a quadric.
\end{corollary}
\begin{proof}
Apply Corollary~\ref{corollary:11-and-1}, 
Lemma~\ref{lemma:no-quadric}, Lemma~\ref{lemma:12-on-quadric} and then 
Lemma~\ref{lemma:contains-quadric}.
\end{proof}

\begin{proof}[Proof of Theorem~\ref{maintheorem}]
Due to Corollary~\ref{corollary:10-points} we may assume that 
$11\le\nolinebreak s\le\nolinebreak 12$. Let $s=11$.
If the points $p_1, \ldots, p_{11}$ span a three-dimensional subspace,
they impose independent conditions on forms of degree $3$ in $\P^4$ 
by Lemma~\ref{lemma:11-on-quadric} and Lemma~\ref{lemma:no-quadric}.
If they span the whole $\P^4$, the statement follows from 
Theorem~\ref{theorem:Eisenbud} and Lemma~\ref{lemma:6-points}.

Now let $s=12$. If no $11$ of the points $p_1, \ldots, p_{12}$ 
are contained in a three-dimensional subspace, then
$\Q$-factoriality of $X$ is again implied by Theorem~\ref{theorem:Eisenbud}. 
If $11$ of these points are contained in a three-dimensional subspace,
we are done by Corollary~\ref{corollary:12-points}. 
\end{proof}

\section{Some constructions}\label{section:constructions}

\begin{example}[{see, for example,~\cite[Example~3]{Cheltsov-quartic}, 
\cite[Example~1.21]{Cheltsov-points}, 
and~\cite[Example~6]{Mella}}]
\label{example:construction}
Let $X$ be a quartic in $\P^4$ with homogeneous coordinates $(x_0:\ldots:x_4)$
containing an irreducible quadric.
Then $X$ is described by an equation $QQ'-LC=0$, where $\deg(Q)=\deg(Q')=2$,
$\deg(L)=1$, $\deg(C)=3$ (for brevity we'll also denote the corresponding 
hypersurfaces in $\P^4$ by the symbols $L$ and $C$ and the corresponding 
quadric surfaces by the symbols $Q$ and $Q'$). 
A sufficiently general quadric of this kind has $12$ simple double 
singularities and is apparently non-$\Q$-factorial.

Let $y=Q/L=C/Q'$, $y'=Q'/L=C/Q$. A rational map 
$$f:(x_0:\ldots:x_4)\mapsto (x_0:\ldots:x_4:y)$$ 
gives a birational equivalence 
of $X$ and a variety $Y\subset\P^5_{(x_0:\ldots:x_4:y)}$
given by the equations
\begin{gather*}
yL=Q,\\
yQ'=C.
\end{gather*}
Similarly a map $f':(x_0:\ldots:x_4)\mapsto (x_0:\ldots:x_4:y')$ 
transforms $X$ to a variety $Y'\subset\P^5_{(x_0:\ldots:x_4:y')}$
given by the equations
\begin{gather*}
y'L=Q',\\
y'Q=C.
\end{gather*}

The map $f$ (resp., $f'$) is a composition of one of the two small resolutions 
of $X$, i.\,e. the blow-up of the (non-Cartier) divisor $Q$ (resp.,
$Q'$), and a contraction of the strict transform of
the divisor $Q'$ (resp., $Q$) to a singular point on $Y$. 
If $Q'$ is nonsingular, 
then the corresponding point on $Y$ is a simple double point. 
A composition $f'\circ f^{-1}:Y\dasharrow Y'$ is a link of type II 
(see~\cite{Matsuki} for a definition). 
Varieties $Z$ and $Z'$ are connected by a flop $\chi$.

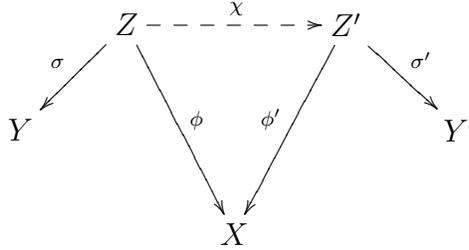
$$
\xymatrix{
&Z\ar@{-->}[rr]^{\chi}\ar[ld]_{\sigma}\ar[rdd]^{\phi}&&
Z'\ar[rd]^{\sigma'}\ar[ldd]_{\phi'}&\\
Y&&&&Y'\\
&&X&&
}
$$

Similarly, if we have a general intersection $Y$ of a quadric $D_2$
and a cubic $D_3$ in $\P^5$
with a unique simple double point $p$, then $Y$ contains $12$ lines
$l_1, \ldots, l_{12}$ passing through $p$, all of them contained 
in the common tangent space $\tilde{L}$ to $D_2$ and $D_3$ at $p$ 
(the generality condition lets us assume that both $D_2$ and $D_3$ are
smooth). Projection $\pi:Y\dasharrow\P^4$ from the point $p$ 
is a composition of a blow-up $\sigma:Z\to Y$ of the point $p$ and 
a contraction $\phi:Z\to X$ of the preimages of the $12$ lines. It is 
easy to see that the image $X$ of the variety $Y$ under this projection
is a quartic, that it contains a quadric 
$\pi(\tilde{L}\cap D_2)$ and that it is singular exactly at the $12$
points $\pi(l_i)$.
\end{example}

\begin{remark}
If the quadric $D_2$ from Example~\ref{example:construction}
degenerates to a cone over a smooth three-dimensional 
quadric with a vertex at $p$, 
and the cubic $D_3$ passes through $p$ and remains sufficiently general,
then $Y$ has a unique simple double point and is birationally equivalent
to a non-$\Q$-factorial double quadric ramified over a quartic section;
the double cover structure is given by the projection $\pi$. Note 
that this gives an example of a non-$\Q$-factorial nodal 
double quadric ramified over a quartic section with the
least possible number of singularities. 
\end{remark}

\begin{remark}[{cf.~\cite[Conjecture~9]{BrCoZu}}]
If the quartic $X$ containing a quadric is sufficiently 
general, then it has no (biregular) automorphisms interchanging the
quadrics $Q$ and $Q'$. 
Indeed, each automorphism of $X$ preserves the elliptic curve 
$E=Q\cap Q'$ and acts on a set 
$\QQ=\{q_1, \ldots, q_{16}\}$ of points $q_i\in E$ such that $4q_i$ is 
a hyperplane section of $E$; if the quadrics $Q$ and $Q'$ are sufficiently 
general, then all the fixed points of the automorphisms from $\Aut(E)$,
acting on $\QQ$, are contained in $\QQ$. Let 
$\tau_1, \ldots, \tau_N\in\PGL_4$ be all automorphisms of $\P^3$ that 
act on the set $\QQ$ and have no fixed points on $E\setminus\QQ$;
for $t\in E\setminus\QQ$ let
$t_j=\tau_j(t)$, $1\le j\le N$. Choose a cubic $C$ so that $C$ 
passes through $t$ but not through any of the points $t_j$.
Assume that there is a nontrivial automorphism $\tau\in\Aut(X)$.
Let $\tau'\in\PGL_4$ be such an automorphism of the subspace $L$ that 
$\left.\tau\right|_L=\left.\tau'\right|_L$. Then $\tau'$ must 
coincide with one of the automorphisms $\tau_j$. Moreover,
it must act on the set of the singular points of $X$, and one of these is $t$. 
Hence $\tau'=id$.

Now let $X$ be a quartic without automorphisms interchanging 
$Q$ and $Q'$. Then the varieties $Y$ and $Y'$ are not isomorphic since
any isomorphisms of $Y$ and $Y'$ gives rise to an automorphism of
$X$ interchanging $Q$ and $Q'$. Hence a general complete intersection 
of a quadric and a cubic in $\P^5$ with one singular point is 
not birationally rigid. On the other hand, one can easily 
construct an example of a quartic $X$ corresponding to isomorphic $Y$ and $Y'$ 
(see~\cite[Example~6]{Mella}).
\end{remark}

A complete intersection of a quadric and a cubic in $\P^5$ is a popular
object of study. A general variety $\Y$ of this type 
is proved to be birationally rigid, and the group of its birational 
automorphisms $\Bir(\Y)$ 
is completely described (see~\cite{Iskovskikh-rigid} 
and~\cite{IskovskikhPukhlikov}).
This group is a semidirect product of its subgroup of biregular automorphisms
$\Aut(\Y)$ and a subgroup generated by birational involutions corresponding to
lines and conics with linear spans contained in the quadric. 
The former involutions come from the double cover structures,
and the latter correspond to elliptic fibrations. 
The following example shows that there are some more birational involutions
in the singular case.

\begin{example}\label{example:involution}
Let $B\subset Y$ be a conic passing through the singular point $p$
such that its linear span is not contained in $D_2$ (a general $Y$ contains
$240$ conics of this type --- see Remark~\ref{remark:identification}).
A general three-dimensional subspace containing $B$ intersects $D_2$ by 
a nonsingular quadric $S$, and $B$ is a curve of bidegree $(1, 1)$ on $S$. 
The intersection of $D_3$ with $S$ consists of $B$ and an elliptic curve
$B'$ of bidegree $(2, 2)$, passing through $p$ and intersecting $B$ 
by three more points. Let $g:\nolinebreak\tilde{Y}\to\nolinebreak Y$ 
be a composition of blow-ups
of $p$ and of a strict transform of the conic $B$; let $E$ and $F$ 
be the corresponding exceptional divisors. Then $\tilde{Y}$ is endowed by
a structure of an elliptic fibration 
$\phi:\tilde{Y}\to\P^2$ with a section $E$ and a $3$-section $F$. 
A reflection with respect to the section $E$ (defined at least on smooth 
fibers) gives rise to a birational automorphism $\tilde{\tau}_B$ 
of $\tilde{Y}$ and to a birational involution $\tau_B$ of $Y$.
The action of $\tilde{\tau}_B$ on $\Pic(\tilde{Y})$ is described by
Statement~\ref{statement:action}.
\end{example}

\begin{statement}[{cf.~\cite[Lemma~5.1.3]{IskovskikhPukhlikov}}] 
\label{statement:action}
Let $h$, $f$ and $e$ denote the classes of $H$, $F$ and $E$ in 
$\Pic(\tilde{X})$. Then
\begin{gather*}
\tilde{\tau}_B^*h=15h-8f-16e,\\
\tilde{\tau}_B^*f=14h-7f-16e,\\
\tilde{\tau}_B^*e=e.
\end{gather*}
\end{statement}
\begin{proof}
The third equality is obvious. Let $G$ be a general fiber of the fibration 
$\phi$; let $h_G$, $f_G$ and $e_G$ denote the restrictions of $h$, $f$ and $e$ 
to $G$.
Note that $\tilde{\tau}_B$ acts on $\Pic(G)$, and
$$\tilde{\tau}_B^*h_G=8f_G-h_G+m(h_G-f_G-e_G)$$
for some $m\in\Z$ since the kernel of the restriction map 
$\Pic(\tilde{X})\to\Pic(G)$ is generated by the class $h_G-f_G-e_G$.

Let $S=\phi^{-1}(l)$ denote a preimage of a general line $l\subset\P^2$;
let $h_S$, $f_S$ and $e_S$ be restrictions of $h$, $f$ and $e$ to $S$. 
Then
\begin{gather*}
h_S^2=4, h_Sf_S=0, h_Se_S=2,\\
f_S^2=e_S^2=-2, f_Se_S=1.
\end{gather*}
Since $(\tilde{\tau}_B^*h_S)^2=h_S^2=6$, we have
$$6=((8f_S-h_S)+m(h_S-f_S-e_S))^2=
-122+8m,$$
and hence $m=16$ that gives the first equality.

Finally, the second inequality is implied by the other two 
since the involution $\tilde{\tau}_B$ preserves the class $h-e-f$.
\end{proof}

\begin{remark}\label{remark:identification}
If $X$ is sufficiently general, then through any of its singular points 
$p_i$, $i=1, \ldots, 12$, there pass $24$ lines $l_{i, j}$, 
$j=1, \ldots, 24$; two of them (say, $l_{i, 1}$ and $l_{i, 2}$) are 
contained in the quadric $Q'$, two (say, $l_{i, 3}$ and $l_{i, 4}$) 
are contained in the quadric $Q$, and other $20$ are not contained even in the 
subspace $L$. Each of the points $p_i$ corresponds to a birational 
involution $\sigma_i$ of $X$, and each line
$l_{i, j}$ corresponds to a birational involution $\sigma_{i, j}$ 
(see~\cite{Pukhlikov-quartic} or~\cite{IskovskikhPukhlikov}). 
On the other hand, each conic $B\subset Y$ such that its linear span 
is contained in $D_2$ corresponds to a birational involution $\tau_B$ 
of $Y$ (see.~\cite{Iskovskikh-rigid} or~\cite{IskovskikhPukhlikov}); each
conic $B\subset Y$ such that its linear span is not contained in $D_2$, and 
$p\in B$, also corresponds to some birational involution $\tau_B$ 
(see Example~\ref{example:involution}). The images of the lines 
$l_{i,3}$ and $l_{i, 4}$ under the map $f$ are conics 
$B_{i, 3}$, $B_{i, 4}$ with linear spans contained in  
$D_2$, and the images of the lines $l_{i, j}$, $5\le j\le 24$, are 
conics $B_{i, j}$ with linear spans not contained in $D_2$. 
All the conics $B_{i, j}$, $1\le i\le 12$, $3\le j\le 24$, 
pass through the singular point $p$. 
Finally, let $l_i\subset Y$ be a line such that 
$\phi(l_i)=p_i$; each line $l\subset Y$ corresponds to 
a birational involution $\tau_l$ (see~\cite{Iskovskikh-rigid} 
or~\cite{IskovskikhPukhlikov}).
Under the natural isomorphism of the groups 
$\Bir(X)$ and $\Bir(Y)$ the involutions $\sigma_i$ are identified
with the involutions $\tau_{l_i}$, $1\le i\le 12$,
and the involutions $\sigma_{i, j}$ are identified with the involutions 
$\tau_{B_{i, j}}$, $1\le i\le 12$, $3\le j\le 24$.
\end{remark}

%%%%%%%%%%%%%%%%%%%%%%%%%%%%%%%%%%%%%%%%%%%%%%%%%%%%%%%%%%%%%%%%%%%%

\end{document}